\documentclass[11pt, twoside, fleqno, draft]{article} 

\usepackage{stmaryrd}
\usepackage{amsmath}
\usepackage{a4}
\usepackage[francais, english]{babel}
\usepackage{euscript}
\usepackage{amsfonts}
\usepackage{amssymb}
\usepackage{xypic}
\usepackage[dvips]{graphics}
\usepackage{enumerate}
\usepackage{array}
\usepackage{theorem}
\usepackage{epic}
\usepackage{color}
\xyoption{all}
\usepackage[T1]{fontenc}
\usepackage{makeidx}
\usepackage{float}

\DeclareMathAlphabet{\mathpzc}{OT1}{pzc}{m}{it}

\newtheorem{theorem}{Th\'eor\`eme}[section]
\newtheorem{proposition}[theorem]{Proposition}
\newtheorem{corollary}[theorem]{Corollaire}
\newtheorem{lemma}[theorem]{Lemme}

\theorembodyfont{\upshape}
\newtheorem{definition}[theorem]{D\'efinition}
\newtheorem{remark}[theorem]{Remarque}

\if not
\renewcommand*{\thesection}{%
    \ifnum\value{section}>0 \arabic{section}\fi%
  }%
\renewcommand*{\thesubsection}{%
    \thesection%
    \ifnum\value{subsection}>0 .\arabic{subsection}\fi%
  }%

\fi

\numberwithin{equation}{section}

\newcommand{\qed}{\hfill $\Box$ \medskip}

\newcommand{\PP}{{\mathbb P}}
\renewcommand{\AA}{{\mathbb A}}

\renewcommand{\O}{{\cal O}}

\newcommand{\Spec}{{\operatorname{Spec}\kern 1pt}}
\newcommand{\Spf}{{\operatorname{Spf}\kern 1pt}}

\def\buildrel#1\over#2{\mathrel{\mathop{\kern 0pt#2}\limits^#1}}

\newcommand{\Ext}{{\mathrm{Ext}}}
\newcommand{\ext}{{\mathpzc{Ext}}}
\newcommand{\Hom}{{\mathrm{Hom}}}

\renewcommand{\H}{{\mathrm{H}}}
\newcommand{\R}{{\mathrm{R}}}
\renewcommand{\Im}{\operatorname{Im}\kern 1pt}
\newcommand{\coker}{\operatorname{coker}\kern 1pt}

\renewcommand{\epsilon}{\varepsilon}
\newcommand{\Aut}{{\mathrm{Aut}}}

\newcommand{\preuve}{\noindent{\it Preuve : }}

\newcommand{\Ar}{{\mathbf{Ar}\kern 0.5pt}}

\newcommand{\cd}{\operatorname{cd}\kern 1pt}

\renewcommand{\tilde}{\widetilde}

\newcommand{\U}{{\mathcal U}}

\newcommand{\p}{{\mathfrak p}}

\newcommand{\A}{{\mathfrak A}}

\newcommand{\m}{{\mathfrak m}}

\newcommand{\F}{{\mathcal F}}
\newcommand{\G}{{\mathcal G}}

\newcommand{\dif}{{\mathfrak d}}

\def\d#1{{\partial \over \partial #1}}

\newcommand{\cart}{\ar@{}[dr] |{\square}}
\newcommand{\comm}{\ar@{}[dr] |{\circlearrowleft}}

\begin{document}

\title{Théorie des déformations équivariantes des morphismes
localement d'intersections complètes}

\author{Sylvain Maugeais}

\maketitle

\selectlanguage{english}

\begin{abstract}
  This is an expository paper on the subject of the title. It assumes basic scheme theory, commutative and homological algebra. 
\end{abstract}

\selectlanguage{francais}

\section{Introduction}

Le but de cet article est de donner les propriétés basiques et les
théorème fondamentaux de la théorie des déformations équivariantes des
schémas d'intersections complètes munies d'une action d'un groupe
abstrait fini. À proprement parler, il n'y a rien de nouveau dans cet
article. Le but essentiel est de donner des preuves simples à des
théorèmes déjà connus et de les illustrer par des exemples. Nous
n'hésiterons donc pas à ajouter des hypothèses sur les objets
considérés, ni à donner des énoncés moins forts, si l'exposé peut y
gagner en simplicité. Des exposés beaucoup plus généraux peuvent être
trouvés dans \cite{IllusieI, IllusieII}, \cite{Laudal} ou
\cite{DeformationWewers}. Les techniques utilisées ici s'inspire de
celles utilisées dans \cite{Vistoli} et \cite{CotangentComplex}.

Nous utilisersons constamment des schémas munis de l'action d'un groupe
fini, il est donc naturel de poser les définitions suivantes.

\begin{definition}
Soient $G$ un groupe fini et $Y$ un schéma. Un \textit{$Y[G]$-schéma} 
(ou $A[G]$-schéma si $Y=\Spec A$) est un couple $(X\to Y, i)$ où $X\to Y$
 est un morphisme de schémas et $i:G\to \Aut_{Y} (X)$ est un homomorphisme de
groupes abstraits. S'il n'y a pas de confusion possible, on notera
aussi $(X\to  Y, G)$, voire $(X, G)$, le couple $(X\to Y, G)$.

Les $Y[G]$-schémas forment une catégorie, les flèches étant les
morphismes $G$-équivariants de $Y$-schémas.

Si $P$ est une propriété des morphismes de schémas, nous dirons qu'un
morphismes de $Y[G]$-schémas $(X_{1}, i_{1})\to (X_{2}, i_{2})$
possède la propriété $P$ si le morphisme sous-jacent $X_{1}\to X_{2}$
possède la propriété $P$.
\end{definition}

Soient $k$ un corps, $G$ un groupe fini et $(X_{0}, i_{0})$ un
$k[G]$-schémas. Une déformation équivariante de $(X_{0}, i_{0})$
consiste en la donnée   
\begin{enumerate}[-]
\item d'un anneau local artinien $A$ d'idéal maximal $\m_A$ de corps
résiduel $k$ ;
\item d'un $A[G]$-schémas $(X,i)$ plat sur $A$ ; 
\item et d'un isomorphisme $(X \times_{\Spec A} \Spec k, i) \to
(X_{0}, i_{0})$
\end{enumerate}
Une telle donnée sera appelée \emph{déformation équivariante} de $(X_{0},
i_{0})$ au-dessus de $A$.

Deux déformations équivariantes $(X_{1}, i_{1})$ et $(X_{2}, i_{2})$
de $(X_{0}, i_{0})$ au-dessus d'un anneau $A$ seront dites isomorphes
s'il existe un $A[G]$-morphisme $f:(X_{1}, i_{1})\to (X_{2}, i_{2})$
induisant l'identité sur les fibres fermées (qui sont canoniquement
isomorphe à $(X_{0}, i_{0})$). On peut remarquer qu'alors $f$ est un
isomorphisme grâce à \cite{Deformation}, Lemma 3.3.

Nous noterons $\A$ la catégorie des anneaux locaux artiniens de corps résiduel $k$.

\begin{definition}[relèvement équivariant]\label{EquivariantLifting}
Soient $A'$ un objet de $\A$, $\m_{A'}$ son idéal maximal,
$\mathfrak{a}\subset A'$ un idéal tel que $\m_{A'}\mathfrak{a}=0$
(en particulier, $\mathfrak{a}$ est naturellement une $k$-algèbre) ; notons 
$A=A'/\mathfrak{a}$ et $\m_A$ son idéal maximal. Soit $(X, i)$ une
déformation équivariante au-dessus de $A$ d'un $k[G]$-schémas $(X_{0},
i_{0})$. Un relèvement équivariant de $(X, i)$ à $A'$ est une
déformation équivariante $(X', i')$ de $(X_{0}, i_{0})$ au-dessus de
$A'$ telle qu'il existe un isomorphisme $(X'\times_{\Spec A'} \Spec A, 
i') \to  (X, i)$.
\end{definition}

L'exemple le plus simple de donnée satisfaisant les hypothèses de la
définition précédente est donné par $A'=k[\varepsilon]/
(\varepsilon^{2})$ (l'algèbre des nombres duaux qu'on note souvent
$k[\varepsilon]$), $\m_{A'}=\mathfrak{a}= (\varepsilon)$ et $A=k$. Les
relèvements équivariants de $(X_{0}, i_{0})$ à $k[\varepsilon]$ seront
appelés les \textit{déformations équivariantes du premiers ordres}.

Le but de cet article est de démontrer une version affaiblie du théorème suivant.

\begin{theorem}\label{Principal}
Soient $k$ un corps, $G$ un groupe fini et $X_{0}$ un $k$-schéma
réduit, noethérien,  localement d'intersection complète et muni d'un
action de $G$ (i.e. il existe un homomorphisme $i_{0}:G\to
\Aut_{k} (X_{0})$ faisant de $(X_{0}, i_{0})$ un $k[G]$-schéma). 
Notons $\Omega_{X_{0}/k}$ le complété séparé du faisceau des différentielles
relatives si $X_{0}$ est le spectre d'un anneau local complet, et notons
$\Omega_{X_{0}/k}$ le faisceau des différentielles relatives
sinon. Reprenons les notations introduites dans la définition 
\ref{EquivariantLifting}. 
Soit $(X, i)$ une déformation équivariante de $(X_{0}, i_{0})$
au-dessus de $A$.
\begin{enumerate}
\item Si $(X', i')$ est un relèvement équivariant de $(X, i)$ à $A'$
alors le groupe des automorphismes de relèvements de  $(X', i')$
(i.e. les automorphismes équivariants de $(X', i')$ qui fixent $(X,
i)$) est canoniquement isomorphe à
$$\mathfrak{a}\otimes_{k} \Hom_{\O_{X_0}} (\Omega_{X_0/k}, \O_{X_0})^G.$$
\item Il existe un élément canonique
$$\omega_{A', \mathfrak{a}}\in\mathfrak{a}\otimes_{k} \Ext^{2}_{G,
\O_{X_0}} (\Omega_{X_0/k}, \O_{X_0})$$ appelé \emph{obstruction},
tel que $\omega_{A', \mathfrak{a}}=0$ si et seulement s'il existe un
relèvement de $(X, i)$ à $A'$. 
\item S'il existe un relèvement de $(X, i)$ à $A'$, alors l'ensemble
des classes d'iso\-morphismes de relèvements de $(X, i)$ à $A'$ est
canoniquement isomorphe à  
$$\mathfrak{a}\otimes_{k} \Ext^{1}_{G,\O_{X_0}} (\Omega_{X_0/k}, \O_{X_0}).$$
\end{enumerate}
\end{theorem}

La définition de morphisme localement d'intersection complète 
adoptée ici (cf. définition \ref{MorphismeGIC}), est un peu plus générale
que celle utilisée habituellement (cf. par exemple \cite{LiuBook})
afin d'englober le cas des anneaux locaux.

Voici le plan de cette note. Dans la section 2, nous rassemblons
les définitions et propriétés concernant les intersections locales
complètes munies d'une action d'un groupe fini. Dans la section 3,
nous montrons le résultat pour les intersections locales complètes
affines (cf. Théorème \ref{IdentificationCasAffine} et Théorème
\ref{ObstructionCasAffine}) et nous démontrons une version plus faible
du cas général dans la section 4 (cf. Théorème \ref{IdentificationCasGeneral}
et Théorème \ref{ExistenceFaibleObstruction}). La section 6 regroupe divers 
exemples et corollaires.  

L'outil de base est la cohomologie équivariante. Notre référence dans
ce domaine est \cite{Tohoku} que nous utiliserons librement.

Le groupe nul sera noté $0$.

\section{Morphismes localement d'intersection complètes}

\begin{definition}\label{MorphismeGFL}
Soient $A$ un anneau et $G$ un groupe fini. Nous dirons qu'un
$A[G]$-schéma est $G$-FL s'il est d'une des formes suivantes
\begin{enumerate}[-]
\item $(\AA^{n}_{A, G}, G)$ où $\AA^{n}_{A, G}:=\Spec A[\left\{x_{i,
\sigma} \right\}_{1\le i \le n, \sigma \in G}]$ et l'action (à gauche) de $G$ est
définie par $\sigma' (X_{i, \sigma})=X_{i, \sigma \sigma'}$ ; 
\item un schéma affine $(\Spec B, G)$, où $B$ est une localisation de 
$\O_{\AA^{n}_{A, G}}$ stable sous l'action de  $G$ ;
\item $(\Spec A[[\left\{x_{i, \sigma}\right\}_{1
\le i \le n, \sigma \in G}]], G)$ et l'action est définie par $\sigma'
(X_{i, \sigma})=X_{i, \sigma \sigma'}.$ 
\end{enumerate}
\end{definition}

La proposition suivante justifie la dénomination $G$-FL (pour
$G$-formellement lisse). 

\begin{proposition}\label{FormallySmoothMorphism}
Soient $A$ un anneau, $G$ un groupe fini et $(Y, G)$ un $A[G]$-schéma
$G$-FL. Soit $(X_{0}, G)$ et $(X, G)$ des $A[G]$-schémas
affines et supposons qu'il existe un diagramme commutatif de $A[G]$-schémas
$$\xymatrix{
(X_{0}, G) \ar[r]^{f}\ar[d] & (Y, G) \ar[d] \\
(X, G) \ar[r] & \Spec A
}$$
où le morphisme $X_{0} \to X$ est une immersion fermée définie par un
idéal de carré nul. Alors il existe un morphisme de $A[G]$--schémas
$\tilde f:(X, G) \to (Y, G)$ rendant commutatif le diagramme
$$\xymatrix{
(X_{0}, G) \ar[r]^{f}\ar[d] & (Y, G) \\
(X, G) \ar[ur]_{\tilde f}. & 
}$$
\end{proposition}

\preuve Immédiate. \qed

Il eut peut être été plus naturel de définir un morphisme $G$-FL comme
étant un morphisme vérifiant la propriété énoncée dans la proposition
\ref{FormallySmoothMorphism}. Toutefois, l'étude de ces morphismes et
de leurs propriétés nous aurait emmenés au-delà des buts de cet article. 

Soit $X\to \Spec A$ un morphisme de schéma. Si $X$ est le spectre du
complété d'un anneau de type fini en un point, nous noterons
$\Omega_{X/A}$ le séparé complété du module des différentielles. Sinon
nous noterons $\Omega_{X/A}$ son faisceau des différentielles relatives. 

En particulier, pour tout faisceau cohérent $\mathcal{F}$ de type
fini, on a une bijection entre l'ensemble des $A$-dérivations de
$\O_{X}$ à valeurs dans $\mathcal{F}$ et l'ensemble
$\Hom_{\O_{X}}(\Omega_{X/A}, \mathcal{F})$.

\begin{proposition}\label{DifferentialIsProjective}
Soient $G$ un groupe fini, $A$ un anneau et $(Y, G)$ un $A[G]$-schéma
$G$-FL. Alors le faisceau $\Omega_{Y/A}$ est un
$\O_{Y,G}$-module projectif.
\end{proposition}

\preuve Ceci découle de la description de $\Omega_{Y/A}$. \qed

\begin{definition}\label{MorphismeGIC}
Soient $G$ un groupe fini, $A$ un anneau et $(X, G)$ un
$A[G]$-schéma. Nous dirons que $(X, G)$ est $G$-IC s'il existe un
$A[G]$-schéma $(Y,G)$ qui soit $G$-FL, et une immersion régulière $(X, G)\to (Y, G)$.
Un morphisme (non équivariant) sera dit d'intersection complète s'il est $0$-IC.
\end{definition}

\begin{remark}
Cette définition de morphisme d'intersection complète n'est pas
classique. Elle englobe évidement la notion de morphisme
d'intersection complète si celui-ci est en plus de type fini. La
définition adoptée ici est nécessaire car, très souvent, nous
``localiserons'' les obstructions et devrons donc travailler sur des
anneaux locaux (voire complet) qui ne sont pas de type fini sur la base.
\end{remark}

En particulier, on voit que si un $A[G]$-schéma $(X, G)$ est $G$-IC et
de présentation finie, alors le schéma sous-jacent $X$ est une
intersection complète. Cette propriété est en fait une équivalence
d'après le lemme suivant.

\begin{lemma}\label{ICImplyGIC}
Soient $A$ un anneau, $G$ un groupe fini et $(X, G)$ un $A[G]$-schéma
affine. Si $X$ est d'intersection complète, alors
$(X,G)$  est $G$-IC.
\end{lemma}

\preuve Supposons que $X\to \Spec A$ est de type fini, la
démonstration dans le cas général se faisant de la même
manière. Notons $B$ l'anneau des fonctions de $X$. Comme $X\to  \Spec
A$ est de type fini, il existe un morphisme 
surjectif 
$$\varphi:C:=A[X_1, \dots , X_{n}]\to B$$
tel que le morphisme induit $X\to \Spec C$ soit une immersion régulière
(cf. \cite{LiuBook}, Corollary 6.3.22). Définissons un morphisme
équivariant
$$\varphi':C':=A[\left\{X_{i, \sigma} \right\}_{1\le i \le n, \sigma \in G}]\to B$$ 
par $\varphi'(X_{i, \sigma})=\sigma (\varphi (X_{i}))$. Comme $\Spec C'\to
\Spec A$ est lisse, le morphisme $\varphi':X\to \Spec C'$ est une
immersion régulière (cf. \cite{LiuBook}, corollary 6.3.22), ce qui achève la preuve. \qed

\begin{lemma}\label{ExactSequence}
Soient $k$ un corps, $(X_{0}, G)$ un $k[G]$-schéma réduit qui soit $G$-IC et
$j:(X_{0}, G)\to (Y_{0}, G)$ une immersion régulière dans un $k[G]$-schéma
qui soit $G$-FL. Notons $\mathcal{I}_{0}$ l'idéal de l'immersion
$j$. Alors on a une suite exacte 
$$0 \to  \mathcal{I}_{0}/\mathcal{I}_{0}^{2} \to
\Omega_{Y_{0}/k}|_{X_{0}}\to \Omega_{X_{0}/k}\to 0.$$
\end{lemma}

\preuve Supposons dans un premier temps que $X_{0}$ et $Y_{0}$ sont de
type fini sur $k$. Alors tout est classique, sauf peut être
l'injectivité. Celle-ci est prouvée dans \cite{Vistoli}, Lemma 4.7
(c'est ici que le fait que $X$ soit réduit intervient).  

La démonstration dans le cas général est similaire.
\qed  

En appliquant le foncteur
$\mathfrak{a}\otimes_{k}\Hom_{\O_{X_0}}(-,\O_{X_0})^G$ à la suite
exacte du lemme précédent on obtient une suite exacte
\begin{multline*}
0\to \mathfrak{a}\otimes_{k} \Hom_{\O_{X_0}}(\Omega_{X_{0}/k},\O_{X_0}) 
\to \mathfrak{a}\otimes_{k} \Hom_{\O_{X_0}}(\Omega_{Y_{0}/k}|_{X_{0}},\O_{X_0})\to  \\
\mathfrak{a}\otimes_{k}\Hom_{\O_{X_0}}(\mathcal{I}_{0}/\mathcal{I}_{0}^{2},\O_{X_0})
\to \mathfrak{a}\otimes_{k} \Ext^{1}_{\O_{X_0}, G}(\Omega_{X_{0}/k},\O_{X_0}) 
\to \\
\mathfrak{a}\otimes_{k} \Ext^{1}_{\O_{X_0},
G}(\Omega_{Y_{0}/k}|_{X_{0}},\O_{X_0}) \to \mathfrak{a}\otimes_{k} \Ext^{1}_{\O_{X_0},
G}(\mathcal{I}_{0}/\mathcal{I}_{0}^{2},\O_{X_0}) \to \\
\mathfrak{a}\otimes_{k} \Ext^{2}_{\O_{X_0}, G}(\Omega_{X_{0}/k},\O_{X_0}) 
\to \mathfrak{a}\otimes_{k} \Ext^{2}_{\O_{X_0}, G}(\Omega_{Y_{0}/k}|_{X_{0}},\O_{X_0})
\end{multline*}

Comme le faisceau $\Omega_{Y_{0}/k}$ est un faisceau de
$\O_{Y_0}[G]$-module projectif (cf. Proposition \ref{DifferentialIsProjective}),
le faisceau $\Omega_{Y_{0}/k}|_{X_{0}}$ est un faisceau de
$\O_{X_0}[G]$-module projectif. Ainsi, on a une suite exacte

\begin{multline}\label{SuiteExacteFondamentale}
0\to \mathfrak{a}\otimes_{k} \Hom_{\O_{X_0}}(\Omega_{X_{0}/k},\O_{X_0}) 
\to \mathfrak{a}\otimes_{k} \Hom_{\O_{X_0}}(\Omega_{Y_{0}/k}|_{X_{0}},\O_{X_0})\to  \\
\mathfrak{a}\otimes_{k}\Hom_{\O_{X_0}}(\mathcal{I}_{0}/\mathcal{I}_{0}^{2},\O_{X_0})
\to \mathfrak{a}\otimes_{k} \Ext^{1}_{\O_{X_0}, G}(\Omega_{X_{0}/k},\O_{X_0}) 
\to 0
\end{multline}

et un isomorphisme
\begin{equation}\label{IsomorphismeFondamental}
\mathfrak{a}\otimes_{k} \Ext^{1}_{\O_{X_0},
G}(\mathcal{I}_{0}/\mathcal{I}_{0}^{2},\O_{X_0}) \to
\mathfrak{a}\otimes_{k} \Ext^{2}_{\O_{X_0}, G}(\Omega_{X_{0}/k},\O_{X_0}) 
\end{equation}

\section{Déformation équivariante des schémas affines}

Fixons quelques notations. Dans cette section et dans la suite $k$
sera un corps, $G$ un groupe fini, $(X_{0}, G)$ un $k[G]$-schéma qui
est $G$-IC, $(X_{0},G)\to  (Y_{0},G)$ une immersion régulière dans un
$k[G]$-schéma qui est $G$-FL. On choisit un objet $A'$ de $\A$,
$\m_{A'}$ son idéal maximal, $\mathfrak{a}\subset A'$ un idéal tel que 
$\mathfrak{a}\m_{A'}=0$. Notons $A=A'/\mathfrak{a}$ et $\m_A$ son
idéal maximal. Finalement, fixons une déformation équivariante $(X,
G)$ de $(X_{0}, G)$ au-dessus de $A$, un $k[G]$-schéma $(Y, G)$
au-dessus de $A$ ayant pour fibre spécial $(Y_{0}, G)$ (on montre
aisément qu'il en existe étant donné la définition des schémas $G$-FL)
et une immersion régulière $\pi:(X, G)\to (Y, G)$ relevant l'immersion
$(X_{0},G)\to (Y_{0}, G)$ (il existe de tels relèvement grâce à la proposition
\ref{FormallySmoothMorphism}). Nous noterons $\mathcal{I}_{0}$ l'idéal
de l'immersion $X_{0}\to Y_{0}$ et $\mathcal{I}$ l'idéal de
l'immersion $X\to Y$. Fixons également un relèvement $G$-FL
$(Y', G)$ de $(Y, G)$ à $A'$.

\subsection{Classification des relèvements}

Soient $(X'_{1}, G)$ et $(X'_{2}, G)$ des relèvements de $(X, G)$ à
$A'$. La proposition \ref{FormallySmoothMorphism} nous fournit des
relèvements $\pi_i:(X'_i, G)\to (Y', G)$ du morphisme $\pi:(X,G)\to(Y,G)$.
On a donc un diagramme commutatif
$$\xymatrix{
\O_{Y'}\ar[rr]^{\pi_1}\ar[rd]\ar[dd] & & \O_{X'_1}\ar[rd]\\
& \O_{Y'}/{\mathfrak{a}}=\O_Y\ar[rr] & & \O_{X} \\
\O_{Y'}\ar[rr]^{\pi_2}\ar[ru] & & \O_{X'_2}\ar[ru]\\
}$$
Notons $\mathcal{I}_{1}$, $\mathcal{I}_{2}$, $\mathcal{I}$ les idéaux de
définitions des immersions régulières $\pi_1$, $\pi_{2}$ et $\pi$.
Comme les $X'_{i}$ sont plats sur $A'$, on a $\mathcal{I}_{i}\otimes_{A'}A=\mathcal{I}$.

Soient $f\in \mathcal{I}$ et $f_{i}\in\mathcal{I}_{i}$ des relèvements
de $f$. Comme $(X'_{1}, G)$ et $(X'_{2}, G)$ sont des relèvements de $(X, G)$ à $A'$
et que $\O_{Y'}$ est plat sur $A'$ (ce qui implique que le noyau de
$\O_{Y'}\to \O_{Y}$ est $\mathfrak{a}\otimes_{A'}\O_{Y'}$) on a
$f_{1}-f_{2}\in \mathfrak{a}\otimes_{A'}\O_{Y'}$.
Or, comme $\m_{A'}\mathfrak{a}=0$, on a un isomorphisme canonique
$$\mathfrak{a}\otimes_{A'} \O_{Y'}\to\mathfrak{a}\otimes_{k}\O_{Y_0}.$$
Ainsi, on peut voir $f_{1}-f_{2}$ comme un élément de
$\mathfrak{a}\otimes_{k}\O_{Y_0}$ qui est nul si et seulement si
$f_{1}=f_{2}$ dans $\O_{Y'}$. 
Notons $\nu_{Y'}((X'_{1}, G), (X'_{2}, G)) (f)$ l'image de $f_{1}-f_{2}$ dans 
$$\mathfrak{a}\otimes_{k}\O_{Y_0}/\mathfrak{a}\otimes_{k}\mathcal{I}_{0}=
\mathfrak{a}\otimes_{k}\O_{X_0}.$$
  
On voit alors aisément que cet élément ne dépend que de $f$, $\pi_{1}$ et
$\pi_{2}$ mais en aucun cas des relèvements $f_{i}$ de $f$.

D'autre part, comme $\pi_{1}$ et $\pi_{2}$ sont équivariants, pour tout
$\sigma\in G$ les éléments $\sigma (f_{i})\in\mathcal{I}_{i}$ sont des
relèvements de $\sigma(f)$. On a donc un morphisme $G$-équivariant de $\O_{Y}$-modules
$$\nu_{Y'} ((X'_{1}, G), (X'_{2}, G)):\mathcal{I}\to\mathfrak{a}\otimes_{k}\O_{X_0}.$$

Finalement, l'isomorphisme canonique 
$$\Hom_{\O_Y} (\mathcal{I}, \mathfrak{a}\otimes_{k}\O_{X_0})^G\cong
\mathfrak{a}\otimes_{k}\Hom_{\O_{X_0}}(\mathcal{I}_{0}/\mathcal{I}_{0}^{2},
\O_{X_0})^{G}$$
(provenant encore une fois de l'égalité $\m_{A'}\mathfrak{a}=0$) nous
permet de considérer $\nu_{Y'} ((X'_{1}, G), (X'_{2}, G))$ comme un
élément de 
$$\mathfrak{a}\otimes_{k}\Hom_{\O_{X_0}}(\mathcal{I}_{0}/\mathcal{I}_{0}^{2},\O_{X_0})^{G}.$$

\begin{proposition}
\label{nuEspacePrincipal}
Soient $(X'_{i}, G)$ ($1 \le i \le 3$) des relèvements de $(X, G)$ à
$A'$. Alors l'application $\nu_{Y'}$ vérifie les propriétés suivantes
(les égalités de schémas doivent être comprises comme des égalités en
tant que sous-schémas fermés de $Y'$)
\begin{enumerate}[{\rm (i)}]
\item $\nu_{Y'} ((X'_{1}, G), (X'_{2}, G))=0$ si et seulement si $(X'_{1},
G)= (X'_{2}, G)$ ; 
\item $\nu_{Y'} ((X'_{1}, G), (X'_{2}, G))+\nu_{Y'} ((X'_{2}, G),
(X'_{3}, G))=\nu_{Y'} ((X'_{1}, G), (X'_{3}, G))$ ; 
\item $\nu_{Y'} ((X'_{1}, G), (X'_{2}, G))=-\nu_{Y'} ((X'_{2}, G),
(X'_{1}, G))$ ; 
\item \label{Transitivite} soit
$\nu\in\mathfrak{a}\otimes_{k}\Hom_{\O_{X_0}}(\mathcal{I}_{0}/\mathcal{I}_{0}^{2},
\O_{X_0})^{G}$, alors il existe un relèvement $(\tilde X_{1}', G)$ de $(X,
G)$ à $A'$ tel que $\nu_{Y'} ((X'_{1}, G), (\tilde X'_{1}, G))=\nu$; 
\item soient $(\tilde Y', G)\to (Y', G)$ un morphisme de
$A'[G]$-schémas tel que $(\tilde Y', G)$ soit $G$-FL, et 
$\tilde \pi_{0}:(X_{0}, G)\to (\tilde Y'\times_{\Spec A'} \Spec k, G)$ une
immersion régu\-lière qui se factorise à travers $(X_{0}, G)\to (Y_{0},
G)$ ; notons $\tilde {\mathcal{I}}_{0}$ l'idéal de l'immer\-sion $\tilde
\pi_{0}$. Alors l'image de $\nu_{\tilde Y'} ((X'_{1}, G), (X'_{2}, G))$
par l'application canonique 
$$\mathfrak{a}\otimes_{k}\Hom_{\O_{X_0}}(\tilde{\mathcal{I}}_{0}/\tilde{\mathcal{I}}_{0}^{2},\O_{X_0})^{G}\to
\mathfrak{a}\otimes_{k}\Hom_{\O_{X_0}}({\mathcal{I}}_{0}/{\mathcal{I}}_{0}^{2},\O_{X_0})^{G}$$ 
est $\nu_{Y'} ((X'_{1}, G), (X'_{2}, G))$.
\end{enumerate}
\end{proposition}

\preuve Seule la partie \eqref{Transitivite} nécessite une
démonstration, les autres propriétés étant des conséquences
immédiates de la définition de $\nu$.

Notons ${\mathcal{I}}_{1}$ l'idéal de l'immersion $X'_{1}\to Y'$,
$\mathcal{I}$ l'idéal de l'immersion $X\to Y$ et 
$$\tilde{\mathcal{I}}_{1}:=\left\{\begin{array}{c}
\tilde f \in\O_{Y'} | \textrm{ l'image }
\bar f \textrm{ de } \tilde f \textrm{ dans }
\O_{Y}=\O_{Y'}\otimes_{A'}A'/{\mathfrak a} \textrm{ est dans }
\mathcal{I} \textrm{ et } \\
\exists f \in \mathcal{I}_{1} \textrm{ tel que l'image de } f-\tilde f
\textrm{ dans } \mathfrak{a}\otimes_{k} \O_{X_0} \textrm{ est } \nu
(\bar f)
\end{array}\right\}$$
On vérifie sans peine que $\tilde {\mathcal{I}}_{1}$ est un idéal de
$\O_{Y'}$ stable sous l'action de $G$ et que le schéma $(\tilde
X'_{1}, G)= (\Spec \O_{Y'}/\mathcal{I}_1, G)$ est un relèvement de $(X,
G)$ tel que $\nu_{Y'} ((X'_{1}, G), (\tilde X'_{1}, G))=\nu$. \qed

Le corollaire suivant est une conséquence immédiate de la proposition ci-dessus.

\begin{corollary}
Soit $(X, G)$ une déformation équivariante de $(X_{0}, G)$. Supposons
qu'il existe un relèvement $(X', G)$ de $(X, G)$ à $A'$. Alors
l'ensemble
$$\left\{\textrm{Relèvement de $(X, G)$ à $A'$ plongés dans  $Y'$} \right\}$$
est canoniquement un espace principal homogène sous l'action de 
$$\mathfrak{a}\otimes_{k}\Hom_{\O_{X_0}}(\mathcal{I}_{0}/\mathcal{I}_{0}^{2},\O_{X_0})^{G}.$$
\end{corollary}

Supposons qu'il existe un isomorphisme de relèvements
$\varphi:(X'_{1}, G)\to (X'_{2}, G)$. On a alors un diagramme
\textbf{non} commutatif \textit{a priori}
$$\xymatrix{
\O_{Y'}\ar[r]^{\pi_{1}^{\#}}\ar[rd]_{\pi_{2}^{\#}} &
\O_{X'_1} \\
& \O_{X'_{2}}\ar[u]_{\varphi^{\#}}.}$$
Posons $D:=\pi_{1}^{\#}-\varphi^\#\circ\pi_{2}^{\#}.$ La composé de $D$
avec la projection $\O_{X'_1}\to
\O_{X'_{1}}\otimes_{A'}A'/\mathfrak{a}$ est nulle (car $\varphi$
induit l'identité sur $X$). Ainsi $D$ est à valeurs dans
$\mathfrak{a}\otimes_{A'}\O_{X'_{1}}$ car $X'_{1}$ est plat sur
$A'$. L'isomorphisme canonique 
$$\mathfrak{a}\otimes_{A'}\O_{X'_{1}}\to
\mathfrak{a}\otimes_{k}\O_{X_0} $$
(provenant de l'égalité $\mathfrak{a}\m_{A'}=0$) permet de considérer
$D$ comme étant à valeurs dans $\mathfrak{a}\otimes_{k}\O_{X_0}$.

Comme les structures de $\O_{Y'}$-modules sur
$\mathfrak{a}\otimes_{k}\O_{X_{0}}$ induites par $\pi_{1}^{\#}$ et
$\varphi^{\#}\circ\pi_{2}^{\#}$ sont les mêmes et que pour tout $f,
g\in\O_{Y'}$ on a 
$$D (fg)=\pi^{\#}_{1} (f) D (g)+ (\varphi^{\#}\circ\pi_{2}^{\#} (g)) D (f),$$
on peut considérer $D$ comme une  $A'$-dérivation de $\O_{Y'}$ à valeurs dans
$\mathfrak{a}\otimes_{k}\O_{X_0}$. Par suite, $D$ induit un morphisme
$$\mu_{Y'} (\varphi):\Omega_{Y'/A'}\to \mathfrak{a}\otimes_{k}\O_{X_0}.$$
D'autre part, comme $\pi_{1}$, $\pi_{2}$ et $\varphi$ sont équivariants,
$\mu_{Y'} (\varphi)$ est en fait un élément de $\Hom_{\O_{Y'}}
(\Omega_{Y'/A'},\mathfrak{a}\otimes_{k}\O_{X_0})^{G}$. Finalement,
l'isomorphisme canonique
$$\Hom_{\O_{Y'}}(\Omega_{Y'/A'},\mathfrak{a}\otimes_{k}\O_{X_0})^{G}\cong
\mathfrak{a}\otimes_{k}\Hom_{\O_{X_{0}}}(\Omega_{Y'/A'}|_{X_{0}},\O_{X_0})^{G}$$
permet de voir $\mu_{Y'} (\varphi)$ comme un élément de
$\mathfrak{a}\otimes_{k}\Hom_{\O_{X_{0}}}(\Omega_{Y'/A'}|_{X_{0}},\O_{X_0})^{G}$.

La proposition suivante est une conséquence immédiate de la
construction de $\mu_{Y'}$.
\begin{proposition}\label{muEspacePrincipal}
Soient $\varphi_1:(X'_1, G) \to (X'_2, G)$ et $\varphi_2:(X'_2, G) \to
(X'_3, G)$ des isomorphismes de relèvements. L'application $\mu_{Y'}$
vérifie les propriétés suivantes  
\begin{enumerate}[{\rm(i)}]
\item $\mu_{Y'}(\varphi_1)=0$ si et seulement si $\varphi_1=Id$.
\item $\mu_{Y'}(\varphi_2 \circ \varphi_1)=\mu_{Y'}(\varphi_2)+\mu_{Y'}(\varphi_1)$.
\item Si $(\tilde Y', G)\to (Y', G)$ est un morphisme entre schémas
$G$-FL et $(X_0, G) \to (\tilde Y_0, G):=(\tilde Y'\times_{\Spec A'}
\Spec k, G)$ est une immersion régulière factorisant le morphisme
$(X_0, G) \to (Y_0, G)$. Alors l'image de $\mu_{\tilde Y'}(\varphi_1)$
par l'application 
$$\mathfrak {a} \otimes_k\Hom_{\O_{X_0}}(\Omega_{\tilde Y_0}|_{X_0},
\O_{X_0})^G \to \mathfrak {a}
\otimes_k\Hom_{\O_{X_0}}(\Omega_{Y_0}|_{X_0}, \O_{X_0})^G$$ est
$\mu_{Y'}(\varphi_1)$.
\item $\mu_{Y'}$ induit une bijection entre les isomorphismes de
relèvements $(X'_1, G)\to (X'_2, G)$ et les préimages de
$\nu_{Y'}((X'_1, G), (X'_2, G))$ par l'application $$\mathfrak {a}
\otimes_k\Hom_{\O_{X_0}}(\Omega_{Y_0/k}|_{X_0}, \O_{X_0})^G \to
\mathfrak {a} \otimes_k\Hom_{\O_{X_0}}(\mathcal I_0/\mathcal I_0^2,
\O_{X_0})^G.$$ 
\end{enumerate}
\end{proposition}

\preuve Le seul énoncé qui n'est pas une conséquence directe de la
définition est le dernier point. Décrivons tout d''abord le morphisme 
$$\mathfrak {a} \otimes_k\Hom_{\O_{X_0}}(\Omega_{Y_0/k}|_{X_0},
\O_{X_0})^G \to \mathfrak {a} \otimes_k\Hom_{\O_{X_0}}(\mathcal
I_0/\mathcal I_0^2, \O_{X_0})^G.$$
En reprenant les notations précédemment introduites, on voit que $D$,
en restriction à $\mathcal{I}$, induit $\nu_{Y'} ((X'_{1}, G),
(X'_{2}, G))$ \textit{via} le morphisme ci-dessus. Il suffit donc de
voir qu'un préimage de $\nu_{Y'} ((X'_{1}, G), (X'_{2}, G))$ induit un
isomorphisme $\varphi:X'_{1}\to X'_{2}$. Fixons donc un tel préimage
$D$. En utilisant les constructions inverses de celles faites
précédement, on peut considérer $D$ comme un morphisme 
$$\O_{Y'}\to \mathfrak{a}\otimes_{A'}\O_{X'_{1}}\to \O_{X'_{1}}.$$
On définit alors un isomorphisme de la manière suivante. Définissons
un morphisme $\tilde \varphi:\O_{Y'} \to \O_{X'_1}$ par
$\tilde\varphi (f)=\pi^{\#}_1(f)+D (f)$. C'est bien un
morphisme d'anneau car $\mathfrak{a}^{2}=0$. Comme $D$ est envoyée sur
$\nu_{Y'} ((X'_{1}, G), (X'_{2}, G))$, on voit que
$\ker\pi_{2}^{\#}\subset \ker \tilde\varphi$  et donc $\tilde\varphi$
induit un morphisme $(X'_{1}, G)\to (X'_{2}, G)$ qui est un isomorphisme par
platitude. \qed

Les propositions \ref{nuEspacePrincipal}, \ref{muEspacePrincipal} et la suite exacte \eqref{SuiteExacteFondamentale} ont pour conséquence le théorème suivant.
\begin{theorem}
\label{IdentificationCasAffine}
Soient $(X_0, G)$ un schéma affine dans la classe $G$-IC, $(X\to\Spec A, G)$ une déformation équivariante de $(X_0, G)$ et $(X'\to\Spec A', G)$ un relèvement de $(X, G)$ à $A'$. Alors 
\begin{enumerate}[{\rm(i)}]
\item le groupe des automorphismes de relèvements de $(X', G)$ est canoniquement isomorphe à $\mathfrak {a} \otimes_k \Hom_{\O_{X_0}}(\Omega_{X_0/k}, \O_{X_0})^G$.
\item l'ensemble des classes d'isomorphie de relèvements de $(X, G)$ à $A'$ est un espace principal homogène sous $\Ext^1_G(\Omega_{X_0/k}, \O_{X_0})$.
\end{enumerate}
\end{theorem}

\subsection{Existence de relèvement}

Fixons un relèvement (\textit{a priori} non équivariant) $X'$ de $X$
(il en existe toujours, cf. par exemple \cite{Vistoli}, Lemma 2.7) et
un plongement $X'\to  Y'$ (ce qui est possible car $Y'$ est
formellement lisse). Notons $\mathcal{I}'$ l'idéal de cet immersion.
Soient $\sigma \in G$ et $f\in \mathcal{I}$. Choisissons un relèvement
$f'$ de $f$ dans $\mathcal{I}'$. Ainsi, $\sigma (f')-f'$ est un
élément de $\O_{Y'}$ (ce n'est pas toujours un élément de $\mathcal{I}'$ car
$X'$ n'est pas forcément muni d'un action de $G$). Notons
$\omega_{X'}(\sigma)(f)$ l'image de $\sigma (f')-f'$ dans
$\O_{X'}$. Comme $X$ est muni d'une action de $G$ (i.e. $\mathcal{I}$
est fixe sous l'action de $G$), $\omega_{X'}(\sigma)(f)$ est nul
modulo $\mathfrak{a}$. On a donc naturellement 
$\omega_{X'}(\sigma)(f)\in\mathfrak{a}\otimes_{A'}\O_{X'}$ (car $X'$
est plat sur $A'$). On montre alors aisément que
$\omega_{X'}(\sigma)(f)$ ne dépend pas du relèvement $f'$ de $f$.

On a donc construit une application $\O_{Y'}$-linéaire
$$\mathcal{I}\to \mathfrak{a}\otimes_{A'}\O_{X'}$$

Les isomorphismes $\mathfrak{a}\otimes_{A'}\O_{X}\cong
\mathfrak{a}\otimes_{k}\O_{X_{0}}$ et
$$\Hom_{\O_{Y}}(\mathcal{I}, \mathfrak{a}\otimes_{k}\O_{X_0})\overset{\sim}{\to }
\mathfrak{a}\otimes_{k}\Hom_{\O_{X_{0}}}(\mathcal{I}_{0}/\mathcal{I}_{0}^2,
\O_{X_0})$$
nous permettent de considérer $\omega_{X'}(\sigma)$ comme un élément
de $$\mathfrak{a}\otimes_{k}\Hom_{\O_{X_{0}}}(\mathcal{I}_{0}/\mathcal{I}_{0}^2,
\O_{X_0}).$$

D'autre part, de la définition de $\omega_{X'} (\sigma)$ découle
directement l'égalité
$$\omega_{X'}(\sigma\sigma')=\sigma (\omega_{X'} (\sigma'))+\omega_{X'}
(\sigma)$$ et donc $\omega_{X'}$ définit un élément de  
$$\H^1(G,
\mathfrak{a}\otimes_{k}\Hom_{\O_{X_{0}}}(\mathcal{I}_{0}/\mathcal{I}_{0}^2,\O_{X_0}))=\mathfrak{a}\otimes_{k}\H^1(G,
\Hom_{\O_{X_{0}}}(\mathcal{I}_{0}/\mathcal{I}_{0}^2,\O_{X_0})).$$
Pour tout élément
$\phi\in\mathfrak{a}\otimes_{k}\Hom_{\O_{X_{0}}}(\mathcal{I}_{0}/\mathcal{I}_{0}^2,\O_{X_0})$
nous noterons $\partial\phi$ le morphisme de cobord définit par $\partial\phi=\sigma\phi-\phi$.

\begin{lemma}
L'image de $\omega_{X'}$ dans  $\mathfrak{a}\otimes_{k}\H^1(G,
\Hom_{\O_{X_{0}}}(\mathcal{I}_{0}/\mathcal{I}_{0}^2,\O_{X_0}))$
est indépendant de $X'$. On la notera $\omega_{emb}$
\end{lemma}

\preuve Choisissons deux relèvements $X'_{1}$ et $X'_{2}$ de $X$.  Le
résultat est une conséquence de l´égalité
$$\omega_{X'_{1}} (\sigma) (f)-\omega_{X'_{2}}(\sigma)(f)=\left( \partial \nu_{Y'} ((X'_{1}, 0), ( X'_{2}, 0))\right) (f).$$ \qed

\begin{proposition}
Il existe un relèvement équivariant de $(X, G)$ à $A'$ si et seulement si $\omega_{emb}=0$.
\end{proposition}

\preuve Le sens direct est immédiat car si on a un relèvement
équivariant $(X', G)$ de $(X, G)$ à $A'$, alors $\omega_{X'}=0$.

Montrons la réciproque et supposons pour cela que
$\omega_{emb}=0$. Choisissons un relèvement quelconque $X'$ de $X$. 
Par définition, il existe $$\phi\in
\mathfrak{a}\otimes_{k}\Hom_{\O_{X_{0}}}(\mathcal{I}_{0}/\mathcal{I}_{0}^2,\O_{X_0})$$
tel que $\omega_{X'}=\partial\phi$.
Choisissons un relèvement (\textit{a priori} non équivariant) $X''$ de
$X$ tel que $\nu_{Y'}((X', 0), (X'', 0))=\phi$. Il est alors aisé de
voir que $\omega_{X''}=0$ et, par suite, qu'on peut munir $X''$ d'une
action de $G$ relevant celle de $X$. \qed

On obtient alors directement le théorème suivant.

\begin{theorem}
\label{ObstructionCasAffine}
Il existe un élément canonique 
$$\omega_{(X, G) , A'}\in
\mathfrak{a}\otimes_{k}\Ext^{2}_{G}(\Omega_{X_{0}/k}, \O_{X_0})$$
qui s'annule si et seulement si $(X, G)$ admet un relèvement
équivariant à $A'$.
\end{theorem}

\preuve C'est une conséquence immédiate de l'isomorphisme
\eqref{IsomorphismeFondamental} et du fait qu'il existe un relèvement
équivariant de $(X, G)$ à $A'$ si et seulement s'il existe un
relèvement équivariant plongé dans $(Y', G)$. \qed

\section{Déformation équivariante des schémas localement d'intersections
complètes}

\begin{definition}
Soient $A$ un anneau, $G$ un groupe fini et $(X, G)$ un
$A[G]$-schéma. Alors $(X, G)$ sera dit localement $G$-IC s'il existe
un recouvrement $\{U_{i} \}$  de $X$ par des ouvers affines stables
sous l'action de $G$ tels que $(U_{i}, G)$ soit un $A[G]$-schéma $G$-IC.
\end{definition}

En particulier, si $(X, G)$ est un $A[G]$-schémas localement $G$-IC
alors le quotient $X/G$ existe. Contrairement au cas des schémas
$G$-IC, il peut exister des $A[G]$-schémas localement d'intersections
complètes qui ne sont pas localement $G$-IC.

\subsection{Classification des relèvements et de leurs automorphismes}

Nous avons tout d'abord besoin de rappeler quelques résultats sur les
extensions équivariantes. 

Soient $G$ un groupe fini, $A$ un anneau, $(X, G)$ un $A[G]$-schéma
noethérien, ${\mathcal F}$ et $\mathcal G$ deux $\O_{X}[G]$-faisceaux
quasi-cohérents sur $X$. Notons $\mathfrak {Ext}_G(\mathcal F,
\mathcal G)$ l'ensemble des classes d'isomorphie d'extensions
équivariantes de ${\mathcal F}$ par $\mathcal G$, i.e. des suites
exactes de $\O_X[G]$-modules
$$0 \to \mathcal G \overset{i}{\to} \mathcal E \overset{j}{\to} \mathcal F \to 0.$$

Il est alors aisé de voir (cf. par exemple \cite{HomologicalAlgebra},
Theorem II.2.4
pour le cas non équivariant)  qu'on a une identification canonique
$$\Ext^1_G(\mathcal F, \mathcal G) \to \mathfrak {Ext}_G(\mathcal F, \mathcal G).$$

Soit $\mathcal U=\{ U_\alpha\}_{\alpha \in \mathcal A}$ un
recouvrement de $X$ par des ouverts stables sous l'action de $G$. On
définit $\mathfrak {Ext}_G(\mathcal U, \mathcal F, \mathcal G)$ comme
l'ensemble des familles  
$$(\{(\mathcal E_\alpha, i_\alpha, j_\alpha)\}_{\alpha \in \mathcal
A}, \{\varphi_{\alpha', \alpha}\}_{\alpha, \alpha' \in \mathcal A})$$ 
avec $(\mathcal E_\alpha, i_\alpha, j_\alpha) \in \mathfrak
{Ext}_G(\mathcal F|_{U_\alpha}, \mathcal G|_{\alpha})$ et  
$$\varphi_{\alpha', \alpha}:(\mathcal E_\alpha, i_\alpha,
j_\alpha)|_{U_\alpha \cap U_{\alpha'}} \to (\mathcal E_{\alpha'},
i_{\alpha'}, j_{\alpha'})|_{U_\alpha \cap U_{\alpha'}}$$ 
 des isomorphismes d'extensions équivariantes vérifiant la condition
de cochaîne 
$\varphi_{\alpha'', \alpha'} \circ \varphi_{\alpha',
\alpha}=\varphi_{\alpha'', \alpha}$ pour tous $\alpha, \alpha',
\alpha'' \in \mathcal A$. 

\begin{lemma}
Soit $\U$ un recouvrement de $X$ par des ouverts affines stables sous $G$.
On a une bijection canonique 
$$\mathfrak {Ext}_G(\mathcal U, \mathcal F, \mathcal
G)\overset{\theta}{\to} \mathfrak {Ext}_G(\mathcal F, \mathcal G)$$  
\end{lemma}

\preuve Soit $(\{(\mathcal E_\alpha, i_\alpha, j_\alpha)\}_{\alpha \in
\mathcal A}, \{\varphi_{\alpha', \alpha}\}_{\alpha, \alpha' \in
\mathcal A}) \in \mathfrak {Ext}_G(\mathcal U, \mathcal F, \mathcal
G)$. Alors, par recollement, il existe une unique extension
équivariante $(\mathcal E, i, j)$ et des
isomorphismes $\psi_\alpha:(\mathcal E, i, j)|_{U_\alpha} \to
(\mathcal E_\alpha, i_\alpha, j_\alpha)$ tels que
$\psi_{\alpha'}=\varphi_{\alpha', \alpha} \circ \psi_{\alpha}$. Posons
$\theta((\{(\mathcal E_\alpha, i_\alpha, j_\alpha)\}_{\alpha \in
\mathcal A}, \{\varphi_{\alpha', \alpha}\}_{\alpha, \alpha' \in
\mathcal A}))=(\mathcal E, i, j)$. 

Nous allons construire une application réciproque. Soit $(\mathcal E,
i, j)\in \mathfrak {Ext}_G(\mathcal F, \mathcal G)$. On a donc une
suite exacte 
$$0 \to \mathcal F \overset{i}{\to} \mathcal E \overset{j}{\to} \mathcal G \to 0.$$
Comme les ouverts $U_\alpha$ sont affines et que $\mathcal F$ est
quasi-cohérent on a une suite exacte 
$$0 \to \mathcal F|_{U_\alpha} \overset{i|_{U_\alpha}}{\to} \mathcal E|_{U_\alpha} \overset{j|_{U_\alpha}}\to \mathcal G|_{U_\alpha} \to 0.$$
Ainsi, $(\mathcal E_\alpha, i_\alpha, j_\alpha):=(\mathcal
E|_{U_\alpha}, i|_{U_\alpha}, j|_{U_\alpha}) \in \Ext^1_G(\mathcal
F|_{U_\alpha}, \mathcal G|_{U_\alpha})$. De plus, comme
$(\mathcal E, i, j)$ est un faisceau en extension, il existe des
isomorphismes  
$$\varphi_{\alpha', \alpha}:(\mathcal E, i, j)_\alpha|_{U_\alpha\cap
{U_\alpha'}} \to (\mathcal E, i, j)_{\alpha'}|_{U_\alpha\cap
U_{\alpha'}}.$$ 
On a donc construit un élément $\theta'(\mathcal E) \in \mathfrak
{Ext}_G(\mathcal U, \mathcal F, \mathcal G)$ dont l'image par $\theta$
est $\mathcal E$. Il est alors aisé de voir que $\theta'$ et $\theta$
sont inverses l'une de l'autre. \qed 

\begin{theorem}
\label{IdentificationCasGeneral}
Soient $k$ un corps, $(X_0, G)$ un $\Spec k[G]$-schéma qui est
localement dans la classe $G$-IC. Soit $A'\to A$ un morphisme
surjectif d'objets de $\A$ et de noyau $\mathfrak{a}$ tel que $\m_{A'}\mathfrak{a}=0$ (où
$\m_{A'}$ désigne l'idéal maximal de $A'$). Soit $(X\to\Spec A, G)$
une déformation équivariante de $(X_0, G)$ et $(X'\to\Spec A', G)$ un
relèvement de $(X, G)$ à $\Spec A'$. 
Alors 
\begin{enumerate}[{\rm(i)}]
\item le groupe des isomorphismes de relèvements de $(X', G)$ est
canoniquement isomorphe à $\mathfrak{a} \otimes_k
\Hom_{\O_{X_0}}(\Omega_{X_0/k}, \O_{X_0})^G$. 
\item l'ensemble des classes d'isomorphie de relèvements de $(X, G)$ à
$A'$ est un espace principal homogène sous $\Ext^1_G(\Omega_{X_0/k},
\O_{X_0})$. 
\end{enumerate}
\end{theorem}

\preuve C'est une conséquence immédiate du cas affine (théorème
\ref{IdentificationCasAffine}) et du lemme précédent. \qed

\subsection{L'obstruction}

Dans le cas général, la construction de l'obstruction nécessite
l'utilisation du complexe cotangent équivariant (cf. \cite{IllusieI,
IllusieII} ou \cite{DeformationWewers}), le problème fondamental
étant que le complexe cotangent habituel, en tant qu'objet de la
catégorie dérivée des complexes de $\O_X$-module, n'est pas muni d'une
action de $G$.  

Nous proposons ici un énoncé plus faible mais qui nous semble plus
explicite. Il repose fondamentalement sur la suite spectrale 
$$\mathrm{I}_{2}^{p, q}:=\H^p(X_{0}/G, \ext^{q}_{G} (\Omega_{X_{0}/k}, 
\O_{X_{0}}))\Rightarrow \Ext^{p+q}_{G} (\Omega_{X_{0}/k},  \O_{X_{0}}).$$
Dans beaucoup de cas concrets, il est plus simple à
manipuler que le théorême général.

\begin{theorem}\label{ExistenceFaibleObstruction}
Soient $k$ un corps, $G$ un groupe fini et $(X_{0}, G)$ un $k[G]$
schéma séparé localement dans la classe $G$-IC. Soient  $A'$ un objet de $\A$,
$\m_{A'}$ son idéal maximal, $\mathfrak{a}\subset A'$ un idéal tel que
$\m_{A'}\mathfrak{a}=0$ ; notons $A=A'/\mathfrak{a}$ et $\m_A$ son
idéal maximal. Soit $(X, G)$ une déformation équivariante au-dessus de
$A$ de $(X_{0}, G)$.

\begin{enumerate}
\item Il existe un élément canonique 
$$\omega_1 \in {\mathfrak{a}}\otimes_{k}\H^0(X_{0}/G, \ext^{2}_{G}
(\Omega_{X_{0}/k}, \O_{X_{0}}))$$ dont l'annulation est nécessaire et suffisante
pour qu'il existe, localement sur $X_{0}$, des déformations
équivariantes de $(X_{0}, G)$.
\item Si $\omega_1=0$ alors il existe un élément canonique
$$\omega_2 \in {\mathfrak{a}}\otimes_{k}\H^1(X_{0}/G, \ext^{1}_{G}
(\Omega_{X_{0}/k}, \O_{X_{0}}))$$ dont l'annulation est nécessaire et suffisante
pour qu'il existe un recouvrement $U_{i}$ de $X_{0}$ et des
relèvements équivariants $(X'_{i}, G)$ de $(X\cap U_{i}, G)$ tels que
$(X'_{i}, G)|_{U_{i}\cap U_{j}}\cong (X'_{j}, G)|_{U_{i}\cap U_{j}}$
\item Si $\omega_1=0$ et $\omega_{2}=0$ alors il existe un élément canonique
$$\omega_3\in {\mathfrak{a}}\otimes_{k}\H^2(X_{0}/G, \ext^{0}_{G}(\Omega_{X_{0}/k},
\O_{X_{0}}))$$ dont l'annulation est nécessaire et suffisante pour
qu'il existe un relèvement équivariant de $(X, G)$ à $A'$.
\end{enumerate}
\end{theorem}

\preuve Démontrons tout d'abord le point 1. Soit $\{ U_{i} \}$ un
recourvrement de $X_{0}$ par des ouverts affines stables sous l'action
de $G$. D'après le théorème \ref{ObstructionCasAffine}, il existe un élément canonique
$\omega_i\in {\mathfrak{a}}\otimes_{k}\H^0(U_{i}/G,
\ext^{2}_{G}(\Omega_{X_{0}/k}, \O_{X_{0}}))$ dont la nullité est
nécessaire et suffisante pour qu'il existe un relevement équivariant de $(X\cap
U_{i}, G)$. Comme cet élément est définit de manière canonique, les
différents $\omega_i$ se recolle pour donner un élément $\omega_{1}\in
{\mathfrak{a}}\otimes_{k}\H^0(X_{0}/G, \ext^{2}_{G}
(\Omega_{X_{0}/k}, \O_{X_{0}}))$.

\medskip
Passons à la démonstration du point 2 et supposons que
$\omega_{1}=0$. Fixons un recouvrement $U_{i}$ de $X_{0}$ par des
ouverts affines stables sous l'action de $G$ et d'intersections
complètes sur $k$. Comme $\omega_1=0$, $\omega_{1}|_{U_{i}}=0$ et donc
il existe des relèvements équivariants $(X'_{i}, G)$ de $(X \cap
U_{i}, G)$. Comme $X_{0}$ est séparé, $U_{i} \cap U_{j}$ est affine et
d'intersection complète sur $k$. Notons $\nu_{i, j}$ l'image de
$\nu_{Y'} ((X'_{i}, G), (X'_{j}, G))$ dans $\mathfrak{a}\otimes_{k}\Ext^{1}_{G}
(\Omega_{U_{i}\cap U_{j}/k}, \O_{U_{i}\cap U_{j}})$ (pour un schéma
$(Y', G)$ G-FL sur $A'$ tel qu'il existe une immersion régulière
équivariante $(X\cap U_{i}\cap U_{j}, G)\to (Y', G)\times_{\Spec A'} \Spec A$ ; l'élément
$\nu_{i, j}$ est indépendant du choix de $Y'$ grâce à la proposition
\ref{nuEspacePrincipal} (v)). La
propriété (ii) de la proposition \ref{nuEspacePrincipal} permet de
voir que les $\nu_{i, j}$ définissent un élément $\omega_{2}$ de 
$${\mathfrak{a}}\otimes_{k}\check \H^1(\{U_{i}\}, \ext^{1}_{G} (\Omega_{X_{0}/k}, \O_{X_{0}}))$$
et on montre aisément que cet élément est indépendant des relèvements
$(X'_{i}, G)$ choisis.

D'autre part, comme les ouverts affines sont acycliques pour la
cohomologie de Zarisky des faisceaux cohérents, on a 
$${\mathfrak{a}}\otimes_{k}\check \H^1(\{U_{i}\}, \ext^{1}_{G} (\Omega_{X_{0}/k}, \O_{X_{0}}))={\mathfrak{a}}\otimes_{k}\H^1(X_{0}, \ext^{1}_{G} (\Omega_{X_{0}/k}, \O_{X_{0}})).$$
On a donc bien définit notre élément $\omega_{1}$ comme dans
l'énoncé. La démonstration du point 2 est alors aisée.

\medskip

Nous allons maintenant démontrer le point 3. Supposons que
$\omega_{1}=0$ et $\omega_2=0$. Comme précédement, fixons un
recouvrement $U_{i}$ de $X_{0}$ par des ouverts affines stables sous
l'action de $G$ et d'intersections complètes sur $k$. Comme
$\omega_{1}=0$ et $\omega_{2}=0$, il existe des relèvements $(X'_{i},
G)$ de $(X\cap U_{i}, G)$ tel que $(X'_{i}\cap (U_{i}\cap U_{j}),
G)\cong (X'_{j}\cap (U_{i}\cap U_{j}),G)$.
Fixons des isomorphismes 
$$\varphi_{i, j}:(X'_{i}\cap (U_{i}\cap U_{j}),G)\to  (X'_{j}\cap (U_{i}\cap U_{j}),G).$$
Notons $\varphi_{i, j, \ell}=\varphi_{i, j}\circ \varphi_{\ell,
j}^{-1}\circ \varphi_{\ell, i} \in
{\mathfrak{a}}\otimes_{k}\Ext^{0}_{G} (\Omega_{U_{i}\cap U_{j}\cap
U_{\ell}/k}, \O_{U_{i}\cap U_{j}\cap U_{\ell}})$ (cf. Théo\-rème
\ref{IdentificationCasAffine}). Les $\varphi_{i, j}$ définissent un élément
$$\omega_{3}\in {\mathfrak{a}}\otimes_{k} \check \H^2 (\{U_{i} \},
\ext^{0}_{G} (\Omega_{X/k}, \O_{x})).$$ La proposition
\ref{muEspacePrincipal} (ii) permet de voir que cet élément est
indépendant des isomorphismes $\varphi_{i, j}$ choisis. On montre alors
aisément que $\omega_{3}=0$ si et seulement s'il existe des
isomorphismes $\varphi_{i, j}$ tels que $\varphi_{i, \ell}=\varphi_{i,
j}\circ\varphi_{j, \ell}$ pour tous $i, j, \ell$. D'où le résultat par
recollement des $X'_{i}$ le long des $U_i\cap U_{j}$ grâce aux
$\varphi_{i, j}$. \qed

\if not

Soient $\F$et $\G$ deux faisceaux de $\O_X[G]$-modules. Considérons
l'ensemble
$$\mathfrak{Ext}^{2}_{G} (\F,\G)=\{0\to \G\to {\mathcal{E}_{1}} \to
{\mathcal{E}_{2}}\to  \F\to  0 \}/\sim ,$$
deux $2$-extensions $P$ et $Q$ étant dites équivalentes s'il existe des
extensions $P_{1}, \dots , P_{n}$ tel que

$$\xymatrix{Q \ar@{~>}[r] & P_{1} \ar@{<~}[r] & P_{2}
\ar@{~>}[r] & \ldots \ar@{<~}[r] & P_{n}\ar@{~>}[r] & P}$$

où on note $\xymatrix{P_{1}\ar@{~>}[r] & P_{2}}$ s'il existe un diagramme
commutatif
$$\xymatrix{
P1 : & 0 \ar[r] & \G\ar[r]\ar[d]_{Id}\ar[r] & {\mathcal{E}_{1}}
\ar[r]\ar[d] & {\mathcal{E}_{2}} \ar[r]\ar[d] & \F\ar[r]\ar[d]^{Id} &
0 \\
P2 : & 0 \ar[r] & \G\ar[r]\ar[r] & {\mathcal{E}_{1}}'
\ar[r] & {\mathcal{E}_{2}}' \ar[r] & \F\ar[r] &  0
}$$

\begin{proposition}
On a une identification naturelle de $\mathfrak{Ext}^{2}_{G} (\F,\G)$
avec $\Ext^{2}_{G} (\F,\G)$.
\end{proposition}

\preuve cf. \cite{HomologicalAlgebra}, IV, Theorem 9.1. \qed

\begin{lemma}
Soit $\{U_{i} \}$ un recouvrement de $X$ par des ouverts affines
stables sous l'action de $G$. Soient $\F$ et $\G$ deux
$\O_X[G]$-faisceaux sur $X$. Alors on a une identification naturelle
de $\Ext^{2}_{G} (\F,\G)$ avec l'ensemble
\end{lemma}

\begin{lemma}\label{PlongementEquivariant}
Soient $X\to \Spec A$ un schéma quasi-projectif et $G$ un groupe fini
agissant sur $X\to \Spec A$. Alors il existe une immersion $X\to
\PP^{n}_{A}$ équivariante.
\end{lemma}

\preuve Il suffit de prouver l'existence d'un $G$-faisceau inversible
très amples sur $X$. Or $X\to \Spec A$ est quasi-projectif donc il
existe un faisceau très amples $\mathcal{L}$. On montre alors aisément
que $\displaystyle{\bigotimes_{g\in G}}g^{*}\mathcal{L}$ est un $G$-faisceau très amples
sur $X$. \qed

\begin{remark}
On peut remarquer que c'est moralement la même démonstration que dans
le cas affine.
\end{remark}

\fi

\section{Corollaires et applications}

Nous allons maintenant voir quelques exemples simples de calculs. La
plupart du temps, nous assumerons une hypothèse de
trivialité de la cohomologie des groupes. Plus précisément, nous
nous placerons presque toujours sous l'une des hypothèses suivantes :
\begin{enumerate}[-]
\item le groupe agissant sur les schémas considérés a un ordre premier
à toutes les caractérisiques intervenant dans le problème ;
\item le groupe agit librement.
\end{enumerate}
Le cas général (sauvage) est beaucoup plus complexe. Pour des exemples
nous renvoyons à \cite{BertinMezard} et \cite{CornelissenKato} pour le
cas des courbes lisses, ou bien \cite{Papier1, Papier2} pour le cas
des courbes stables.

\subsection{Calcul de $\Ext^1_G(\Omega_{X/k},\O_X)$}

\begin{proposition}
Soient $k$ un corps, $X\to \Spec k$ un schéma quasi-projectif
et $G$ un groupe d'ordre premier à la caractéristique de $k$ et
agissant sur $X\to \Spec k$. Notons $\pi:X\to X/G=:\Sigma$ le
morphisme quotient et supposons que $\Sigma\to \Spec k$ est
lisse. Alors on a une suite exacte
\begin{multline*}
\ldots \to \H^{i-1}(\Sigma, \Omega_{\Sigma/k}^{\vee})\to \Ext^{i}_{G} (\Omega_{X/\Sigma}, \O_X)\to \Ext^{i}_{G}(\Omega_{X/k}, \O_X) \to \\
\H^i(\Sigma, \Omega_{\Sigma/k}^{\vee})\to \Ext^{i+1}_{G} (\Omega_{X/\Sigma}, \O_X) \to \ldots 
\end{multline*}
\end{proposition}

\preuve Comme $\pi$ est séparable on a une suite exacte
$$0\to \pi^{*}\Omega^{1}_{\Sigma/k}\to \Omega_{X/k}\to
\Omega_{X/\Sigma}\to 0.$$
On a donc une suite exacte longue donnée par la cohomologie
équivariante
\begin{multline*}
\ldots \to \Ext^{i-1}_{G} (\Omega_{\Sigma/k}, \O_X) \to \Ext^{i}_{G} (\Omega_{X/\Sigma}, \O_X)\to \Ext^{i}_{G}
(\Omega_{X/k}, \O_X) \to \\
\Ext^{i}_{G} (\pi^{*}\Omega_{\Sigma/k}, \O_X)\to \Ext^{i+1}_{G} (\Omega_{X/\Sigma}, \O_X) \to \ldots 
\end{multline*}
Il s'agit alors de calculer $\Ext^{i}_{G} (\pi^{*}\Omega_{\Sigma/k}, \O_X)$. Or on
a une suite spectrale $\textrm{I}^{p, q}_{2}:=\H^p(\Sigma, \ext^{q}_{G}
(\pi^{*}\Omega_{\Sigma/k}, \O_C))$ qui converge vers $\Ext^{p+q}_{G}
(\pi^{*}\Omega_{\Sigma/k}, \O_X)$.
D'autre part, comme $|G|$ est inversible dans $k$, on montre aisément,
en utilisant une suite spectrale, que $\ext^{q}_{G} (\pi^{*}\Omega_{\Sigma/k},
\O_C)=\pi_{*}^{G}\ext^{q}(\pi^{*}\Omega_{\Sigma/k}, \O_C)$. Comme
$\Sigma\to \Spec k$ est lisse, $\Omega_{\Sigma/k}$ est localement
libre donc $\ext^{q} (\pi^{*}\Omega_{\Sigma/k},
\O_C)=0$ si $q>0$. On a donc $\Ext^{i}_{G}(\pi^{*}\Omega_{\Sigma/k},
\O_X)=\H^i(\Sigma, \pi_{*}^{G}\ext^{0} (\pi^{*}\Omega_{\Sigma/k}, \O_X))$.
Or par adjonction, on a un isomorphisme canonique $\pi_{*}\ext^{0} (\pi^{*}\Omega_{\Sigma/k},
\O_X)=\ext^{0} (\Omega_{\Sigma/k}, \pi_{*}\O_X)$. Il en résulte que
$$\pi_{*}^{G}\ext^{0} (\pi^{*}\Omega_{\Sigma/k},
\O_X)=\ext^{0} (\Omega_{\Sigma/k}, \O_\Sigma).$$
D'où le résultat. \qed

\begin{proposition}
Soient $C\to \Spec k$ une courbe lisse sur un corps algébriquement
clos et $G$ un groupe d'ordre premier à la caractéristique de $k$ et
agissant fidèlement sur $C\to\Spec k$. Notons $n$ le nombre de points de
ramification dans le morphisme quotient. Alors on a
$$\dim_k \Ext^{1}_{G} (\Omega_{C/k}, \O_C)-\dim_k \Ext^{0}_{G}
(\Omega_{C/k}, \O_C)=\chi(\Omega_{\Sigma/k}^{\vee})+n.$$ 
\end{proposition}

\preuve D'après la proposition ci-dessus, il suffit de montrer que
$$\sum_i (-1)^{i} \dim_k\Ext^{i}_{G} (\Omega_{C/\Sigma},
\O_\Sigma)=n.$$ 
Comme l'ordre de $G$ est inversible dans $k$, on a
$$\ext^{q}_{G}(\Omega_{C/\Sigma},
\O_\Sigma)=\pi_{*}^{G}\ext^{q}(\Omega_{C/\Sigma}, \O_C)$$ pour tout
$q\ge 0$. 
Soient $\p\in C$ un point ramifié, $D_{\p}$ son stabilisateur et
$\dif_{\p}$ la différente en ce point. Notons $R=\O_{C,\p}$ et
$t$ une uniformisante de $R$. On a alors $\Omega_{C/\Sigma, \p}\cong
\left(R/(t^{\dif_\p}) \right)dt$. On prouve alors aisément (par exemple en prenant
une résolution projective) que $\Ext^{i}(\Omega_{C/\Sigma, \p}, R)=0$
si $i\not = 1$ et $\Ext^{1}(\Omega_{C/\Sigma, \p}, R)\cong \left(R/
(t^{\dif_{\p}})\d{t} \right)$. Il s'agit alors de calculer $\dim_k\left( R/
(t^{\dif_{\p}})\d{t} \right)^{D_{\p}}$ dont on voit (par exemple en
complétant $R$ et en linérisant l'action) qu'il est égal à
$\left\lceil {\dif_{\p}\over |D_{\p}|}\right\rceil$. Or
$\dif_{\p}=|D_\p|-1$ car l'ordre de $D_{\p}$ est premier à la
caractéristique de $k$.
La suite spectrale $\textrm{I}_{2}^{p, q}:=\H^p(\Sigma,
\ext^{q}_{G}(\Omega_{C/\Sigma}, \O_\Sigma))$ 
qui converge vers $\Ext^{p+q}_{G} (\Omega_{C/\Sigma},
\O_\Sigma)$ permet aisément de conclure. \qed

\subsection{Calcul de $\Ext^2_G(\Omega_{X/k},\O_X)$}

\begin{proposition}
Soient $X\to \Spec k$ un morphisme localement d'intersection complète
et $G$ un groupe fini d'ordre premier à la caractéristique de $k$ et
agissant sur $X$. Alors $\Ext^{2}_{G} (\Omega_{X/k}, \O_X)=\Ext^{2}
(\Omega_{X/k}, \O_X)^G$.
\end{proposition}

\preuve On a une suite spectrale $\H^p(G, \Ext^{q} (\Omega_{X/k},
\O_X))$ qui converge vers $\Ext^{p+q}_{G} (\Omega_{X/k}, \O_X)$. Comme
l'ordre de $G$ est inversible dans $k$, 
on a $$\H^p(G, \Ext^{q} (\Omega_{X/k}, \O_X))=0$$ si $p \ge 1$. D'où le résultat. \qed

\begin{proposition}\label{ActionLibreOuvertDense}
Soient $X\to \Spec k$ un schéma affine lisse et $G$ un groupe agissant
librement sur $X\to \Spec k$. Alors $\Ext^{\ell}_{G} (\Omega_{X/k},
\O_X)=0$ pour $\ell \ge 1$.
\end{proposition}

\preuve On a une suite spectrale $\H^p(X/G, \ext^{q}_{G} (\Omega_{X/k},
\O_X))$ qui converge vers $\Ext^{p+q}_{G} (\Omega_{X/k}, \O_X)$. Par
suite, comme $X/G$ est affine et que $\ext^{q}_{G} (\Omega_{X/k}, \O_X)$ est
quasi-cohérent on a $\Ext^{\ell}_{G} (\Omega_{X/k}, \O_X)=\H^0(X/G,
\ext^{\ell}_{G} (\Omega_{X/k}, \O_X))$.
D'autre part, si on note $\pi:X\to X/G$ le morphisme quotient, on a
une suite spectrale $\R^p\pi_{*}^{G}\ext^{q} (\Omega_{X/k}, \O_X)$
qui converge vers $\ext^{p+q}_{G} (\Omega_{X/k}, \O_X)$. Comme $X\to
\Spec k$ est lisse, on a $\ext^{q} (\Omega_{X/k}, \O_X)=0$ pour $q\ge
1$. D'autre part, comme $G$ agit librement, le foncteur $\pi_{*}^{G}$
est exact. On a donc $\ext^{\ell}_{G} (\Omega_{X/k}, \O_X)=0$ pour
$\ell\ge 1$. \qed

\begin{corollary}
Soient $X\to \Spec k$ un schéma affine lisse et $G$  un groupe fini
agissant librement sur $X\to \Spec k$. Pour tout anneau local
aritinien $A$ de corps résiduel $k$ il existe une unique déformation
équivariante de $(X, G)$ au-dessus de $A$ (à isomorphisme près).
\end{corollary}

\preuve Il suffit de combiner le théorème \ref{Principal} (où le
théorème \ref{ObstructionCasAffine}) et la
proposition \ref{ActionLibreOuvertDense}. \qed

\bibliographystyle{../hamsalpha}

\providecommand{\bysame}{\leavevmode\hbox to3em{\hrulefill}\thinspace}

\end{document}